\documentclass[12pt]{amsart}
\usepackage[utf8]{inputenc}
\usepackage{amsmath}
\usepackage{amsfonts}
\usepackage{amssymb}
\usepackage{geometry}
\usepackage{hyperref}
\usepackage{listings}
\usepackage{xcolor}
\usepackage{algorithm}
\usepackage{algpseudocode}
\algrenewcommand\alglinenumber[1]{\scriptsize$\triangleright$}

\hypersetup{
    colorlinks=true,
    linkcolor=blue,
    filecolor=magenta,      
    urlcolor=cyan,
}

\usepackage{tcolorbox}  

\definecolor{codegreen}{rgb}{0,0.6,0}
\definecolor{codegray}{rgb}{0.5,0.5,0.5}
\definecolor{codepurple}{rgb}{0.58,0,0.82}
\definecolor{backcolour}{rgb}{0.95,0.95,0.92}

\lstdefinestyle{SAGEMATHstyle}{
    backgroundcolor=\color{backcolour},   
    commentstyle=\color{codegreen},
    keywordstyle=\color{blue},
    numberstyle=\tiny\color{codegray},
    stringstyle=\color{codepurple},
    basicstyle=\ttfamily\footnotesize,
    breakatwhitespace=false,         
    breaklines=true,                 
    captionpos=b,                    
    keepspaces=true,                 
    numbers=none,                     
    numbersep=5pt,                  
    showspaces=false,                
    showstringspaces=false,
    showtabs=false,                  
    tabsize=2
}
\title[On the computation of $\mathrm{Ext}_{\mathcal A}^{k,k+*}(\mathbb{Z}/2,\mathbb{Z}/2)$ for $k \leq 5$]{On the computation of $\mathrm{Ext}_{\mathcal A}^{k,k+*}(\mathbb{Z}/2,\mathbb{Z}/2)$ for $k \leq 5$}

\author{Phuc Vo Dang$^{*}$}
\address{Department of AI, FPT University, Quy Nhon AI Campus\\
An Phu Thinh New Urban Area, Quy Nhon City, Binh Dinh, Vietnam}
\email{dangphuc150488@gmail.com}

\begin{document}

\maketitle

\begin{abstract}
This Note presents a computational algorithm for determining a basis of the cohomology of the mod 2 Steenrod algebra, $\mathrm{Ext}_{\mathcal A}^{k, k+*}(\mathbb{Z}/2, \mathbb{Z}/2)$ for $k \leq 5$, based on the well-known generators and the Adams relations given in Chen's work~\cite{Chen}. The purpose of this algorithm is to verify the hand-computed results for $\mathrm{Ext}_{\mathcal A}^{4, 4+*}(\mathbb{Z}/2, \mathbb{Z}/2)$ presented in our corrected papers~\cite{Phuc1, Phuc2, Phuc3}. Combining our most recent works \cite{Phuc4, Phuc5}, the verification of the domain and the codomain of the fourth algebraic transfer \cite{Singer} in specific degrees can now be completed.

\end{abstract}


\medskip


\section{Introduction}

The object of our study is the cohomology of the mod 2 Steenrod algebra $\mathcal{A}$, denoted by the Ext groups $\mathrm{Ext}_{\mathcal A}^{k,*}(\mathbb{Z}/2, \mathbb{Z}/2)$. These groups are of fundamental importance in algebraic topology as they form the $E_2$-term of the Adams spectral sequence, a primary tool for calculating the 2-primary stable homotopy groups of spheres. The structure of these Ext groups is notoriously complex and has been studied by many authors (see, e.g., \cite{Bruner, Chen, Hung, Hung2, Chon, Ha, Lin, Phuc1, Phuc2, Phuc3, Singer, Sum, Tangora}). Remarkably, the structure of the indecomposable elements in $\mathrm{Ext}_{\mathcal A}^{k,*}(\mathbb{Z}/2, \mathbb{Z}/2)$ for homological degrees $k \leq 5$ was explicitly determined in~\cite{Chen, Lin, Tangora}.

While the indecomposable elements of $\mathrm{Ext}_{\mathcal A}^{k,*}(\mathbb{Z}/2,\mathbb{Z}/2)$ for $k \leq 5$ are known, computing explicit generators for $\mathrm{Ext}_{\mathcal A}^{k,k+n}(\mathbb{Z}/2,\mathbb{Z}/2)$ at generic degrees $n$ (e.g., $n = 2^{s+t+u} + 2^{s+t} + 2^{s} - 3$), as in our previous works and Nguyen Sum's paper \cite{Sum}, remains error-prone due to the vast number of parameter combinations. Indeed, for instance, in Sum's paper~\cite{Sum}, at degree $n = 2^{s+1} - 3$, a computational error was made in the case $s = 5$ (i.e., $n = 61$), and similarly, at degree $n = 2^{s+1} - 2$, an error occurred for $s = 6$ (i.e., $n = 126$). More precisely, the incorrect results presented in Sum's work ~\cite{Sum} were:
\[
\mathrm{Ext}_{\mathcal A}^{4,\,4+61}(\mathbb{Z}/2,\mathbb{Z}/2) = 0, \quad 
\mathrm{Ext}_{\mathcal A}^{4,\,4+126}(\mathbb{Z}/2,\mathbb{Z}/2) = \langle h_0^{2}h_6^{2} \rangle.
\]
These erroneous results were later corrected in the arXiv version:1710.07895. The accurate values, as verified by the output of our algorithm in the subsequent section, are:
$$\mathrm{Ext}_{\mathcal A}^{4,4+61}(\mathbb{Z}/2,\mathbb{Z}/2) = \langle D_3(0) \rangle, \ \ \mathrm{Ext}_{\mathcal A}^{4,4+126}(\mathbb{Z}/2,\mathbb{Z}/2) = \langle h_0^{2}h^{2}_6, D_3(1)\rangle,$$
where $D_3(s)\neq 0$ for all $s\geq 0.$ Similar computational errors concerning $\mathrm{Ext}_{\mathcal A}^{4,\,4+*}(\mathbb{Z}/2,\mathbb{Z}/2)$ in previously published works have been thoroughly corrected in the online versions~\cite{Phuc1, Phuc2, Phuc3}. 

\medskip

For the above reasons, the purpose of this Note is to present an algorithm, implemented in \textsc{SageMath}, for automatically computing the generators of $\mathrm{Ext}_{\mathcal A}^{k,k+n}(\mathbb{Z}/2,\mathbb{Z}/2)$ for $k \leq 5$, as $n$ varies according to given parameters. This algorithm enables us to rigorously and accurately verify the results previously obtained by manual computations.

\medskip

This paper can be regarded as the ``\textbf{final piece}'', together with our most recent works \cite{Phuc4, Phuc5}, to complete the verification of the dimensions and bases for both the domain and codomain of the fourth Singer algebraic transfer \cite{Singer} in specific degrees. After that, based on a certain structural pattern of admissible monomial basis classes in the quotient space $\mathbb{Z}/2 \otimes_{\mathcal{A}} H^{*}(B(\mathbb Z/2)^{k})$ in each positive degree, we are able to generalize the result and confirm Singer's conjecture on the injectivity of this transfer, as stated in our previous works \cite{Phuc1, Phuc2, Phuc3}.

\section{The algebraic framework}

The algorithm's correctness hinges entirely on the algebraic structure of $\mathrm{Ext}_{\mathcal A}^{*,*}(\mathbb Z/2, \mathbb Z/2)$ as described by Chen \cite{Chen}. This structure consists of two key components: the set of generators and the relations between them.

\subsection{Generators}
The algebra $\mathrm{Ext}_{\mathcal A}^{*,*}(\mathbb Z/2, \mathbb Z/2)$ is generated by a set of indecomposable elements. Chen's paper \cite{Chen} summarizes the known generators for $k \le 5$.
\begin{itemize}
    \item For $k=1$, the generators are $h_i$ for $i \ge 0$.
    \item For $k=3$, the indecomposable generators are the family $c_i$ for $i \ge 0$.
    \item For $k=4$, the indecomposable generators are the seven families $d_i, e_i, f_i, g_{i+1}, p_i, D_3(i),$ and $p_i'$ for $i \ge 0$.
    \item For $k=5$, the indecomposable generators as the thirteen families $$P^1h_1, P^1h_2, n_i, x_i, D_1(i), H_1(i), Q_3(i), K_i, J_i, T_i, V_i, V_i', U_i\ \mbox{for $i \ge 0$}.$$
\end{itemize}
Any element in $\mathrm{Ext}_{\mathcal A}^{k,k+*}(\mathbb Z/2, \mathbb Z/2)$ for $k \le 5$ can be expressed as a polynomial of these generators.

\subsection{Relations}
The generators are not free; they are subject to a large number of relations, which define the algebra's multiplicative structure. Specifically,  Theorem~1.2 in~\cite{Chen} provides a complete list of these relations.
\begin{itemize}
    \item \textbf{General relations:} These include $h_i h_{i+1} = 0$, $h_i^3 = h_{i-1}^2 h_{i+1}$, and $h_i h_{i+2}^2 = 0$.
    \item \textbf{Product relations for $k \le 4$:} These include relations such as $h_j c_i = 0$ for specific indices $j$ (e.g., $j = i-1, i, i+2, i+3$).
    \item \textbf{Product relations for $k=5$:} The most complex part is a set of 39 families of relations that govern the products of generators in ranks $k = 1$, $3$, and $4$, which result in elements of rank $k = 5$. These range from simple vanishing relations like $h_{j+4}^2 c_j = 0$ to complex equalities like $h_{j+4}h_{j+1}c_j = h_{j+3}e_j$.
\end{itemize}

\section{The computational algorithm}

\subsection{Overview}
The algorithm operates as follows:
\begin{itemize}
    \item \textbf{Generate Monomials:} For a target bidegree $(k, k+*)$, the algorithm first generates all possible products of the aforementioned generators (monomials) whose degrees sum to $(k, k+*)$. This forms a set of potential basis elements.
    \item \textbf{Construct Vector Space:} This set of monomials is used to define a basis for a vector space $V$ over $\mathbb{Z}/2$. Each monomial corresponds to a unique basis vector.
    \item \textbf{Encode Relations:} Each relation from Chen's paper is translated into a vector in $V$. For instance, a relation like $A = B+C$ is encoded as the vector corresponding to $A+B+C$, which must be zero in the final group.
    \item \textbf{Compute Kernel:} The set of all relation vectors spans a subspace $R \subseteq V$. The desired Ext group is isomorphic to the quotient space $V/R$. The basis for this quotient space is found by computing the kernel (or null space) of the matrix whose rows are the relation vectors. The dimension of this kernel is the dimension of $\mathrm{Ext}_{\mathcal A}^{k,k+*}(\mathbb Z/2, \mathbb Z/2)$, and its basis vectors correspond directly to the basis elements of the group.
\end{itemize}

\subsection{Code implementation details}

\begin{itemize}
    \item \texttt{get\_generator\_degrees()}: This function acts as the algorithm's data bank. It stores the degree formulas for all indecomposable generators for $k \le 5$, directly implementing the definitions from section (1.1) of Chen's paper. 
    
    \item \texttt{calculate\_ext\_basis(k, n)}: This is the main function. It begins by creating the list of ``potential\_gens'' by iterating through all generator types and indices, forming products, and keeping only those whose total degree matches the target $(k, k+n)$. It then initializes a vector space over \texttt{GF(2)} whose dimension equals the number of these potential generators.
    
    \item \texttt{\_get\_relations()}: This is the core of the implementation. It is a direct, line-by-line coding of the relations listed in  \cite[Theorem 1.2]{Chen}.
    \begin{itemize}
        \item The general $h_i$ relations are implemented first. For example, the relation $h_i^3 = h_{i-1}^2 h_{i+1}$ is translated into the vector sum ``\texttt{get\_vec((('h',i), ('h',i), ('h',i))) + get\_vec((('h',i-1), ('h',i-1), ('h',i+1)))}'', which represents the equation $h_i^3 + h_{i-1}^2 h_{i+1} = 0$.
        \item For $k = 5$, the 39 specific relations from  \cite[Theorem 1.2(iv)]{Chen} are explicitly added. For instance, the relation $h_{j+1}e_j = h_j f_j$ is encoded as the vector ``\texttt{get\_vec((('h',j+1),('e',j))) + get\_vec((('h',j),('f',j)))}''.
    \end{itemize}
    
    \item The collected relation vectors form a ``\texttt{relation\_matrix}''. \\
The crucial step ``\texttt{relation\_matrix.transpose().kernel()}'' computes the basis for the null space of the relation matrix's transpose. This is equivalent to finding the basis for the quotient space $V/R$, yielding the final basis for the Ext group. The verbose output then presents this basis and the simplified relations.
\end{itemize}

\section{Representation of the algorithm}

\noindent\rule{\textwidth}{0.4pt}
\begin{center}
\noindent \textbf{Algorithm: Calculation of $\mathrm{Ext}_{\mathcal A}^{k,\, t = k+n}(\mathbb{Z}/2, \mathbb{Z}/2)$ basis}
\end{center}
\noindent\rule{\textwidth}{0.4pt}

\begin{algorithmic}[1]
\Function{CalculateExtBasis}{$k, n$}
    \State $t \leftarrow k+n$ \Comment{$k$ is the homological degree, and $n$ is the stem degree $t-k$}
    \State $\text{GenDegrees} \leftarrow \text{LoadGeneratorDegreeFormulas()}$ \Comment{From \cite[Equation (1.1)]{Chen}}
    \State $\text{PotentialMonomials} \leftarrow \emptyset$
    \State \Comment{Generate all monomials with homological degree $k$ and internal degree $t$}
    \State $\text{Combinations} \leftarrow \text{GenerateGeneratorCombinations}(s)$
    \ForAll{monomial $m$ in $\text{Combinations}$}
        \If{$\sum_{g \in m} \text{GenDegrees}(g) = t$}
            \State $\text{PotentialMonomials} \leftarrow \text{PotentialMonomials} \cup \{\text{CanonicalForm}(m)\}$
        \EndIf
    \EndFor
    \Statex
    \State $\text{UniqueMonomials} \leftarrow \text{SortAndUnique}(\text{PotentialMonomials})$
    \State $N \leftarrow |\text{UniqueMonomials}|$
    \If{$N=0$} \Return $\emptyset$ \EndIf
    \State $V \leftarrow (\mathbb{Z}/2)^N$ \Comment{Vector space over GF(2)}
    \State $\text{MonomialMap} \leftarrow \text{map from UniqueMonomials to basis vectors of } V$
    \Statex
    \State $\text{RelationVectors} \leftarrow \text{GetRelationVectors}(s, \text{UniqueMonomials}, \text{MonomialMap})$
    \State $M_{\text{rel}} \leftarrow \text{Matrix with rows from RelationVectors}$
    \State $\mathcal{K} \leftarrow \text{Kernel}(M_{\text{rel}})$ \Comment{Null space of the relation matrix}
    \Statex
    \State $\text{Basis} \leftarrow \emptyset$
    \ForAll{\text{vector } $v \in$ \text{Basis}($\mathcal{K}$)}
        \State $\text{basis\_monomial} \leftarrow \text{MonomialMap}^{-1}(v)$
        \State $\text{Basis} \leftarrow \text{Basis} \cup \{\text{basis\_monomial}\}$
    \EndFor
    \State \Return Basis
\EndFunction
\end{algorithmic}

{\bf Generation of Relation Vectors}

\begin{algorithmic}[1]
\Function{GetRelationVectors}{$k, \text{Monomials}, \text{Map}$}
    \State $\text{Vectors} \leftarrow \emptyset$
    \If{$k \ge 3$} \Comment{General $h$-relations from \cite[Theorem 1.2(i)]{Chen}}
        \ForAll{$m \in \text{Monomials}$ containing $h_i^3$ for some $i \ge 1$}
            \State $m' \leftarrow m \text{ with } h_i^3 \text{ replaced by } h_{i-1}^2 h_{i+1}$
            \State $\text{Vectors} \leftarrow \text{Vectors} \cup \{ \text{Map}(m) + \text{Map}(m') \}$
        \EndFor
        \ForAll{$m \in \text{Monomials}$ containing $h_i h_{i+1}$}
             \State $\text{Vectors} \leftarrow \text{Vectors} \cup \{ \text{Map}(m) \}$
        \EndFor
        \Comment{... Other general relations ...}
    \EndIf
    \If{$k \ge 4$} \Comment{Relations with $c_i$ from \cite[Theorem 1.2(ii)]{Chen}}
        \ForAll{$m \in \text{Monomials}$ containing $h_j c_i$ where $j \in \{i-1, i, i+2, i+3\}$}
            \State $\text{Vectors} \leftarrow \text{Vectors} \cup \{ \text{Map}(m) \}$
        \EndFor
    \EndIf
    \If{$k = 5$} \Comment{The 39 relations from \cite[Theorem 1.2(iv)]{Chen}}
         \ForAll{indices $j$ in a suitable range}
            \State $m_1 \leftarrow \text{CanonicalForm}(h_{j+4}^2 c_j)$ \Comment{Example relation: $h_{j+4}^2 c_j = 0$}
            \If{$m_1 \in \text{Map}$} $\text{Vectors} \leftarrow \text{Vectors} \cup \{ \text{Map}(m_1) \}$ \EndIf
            
            \State $m_A \leftarrow \text{CanonicalForm}(h_{j+4}h_{j+1}c_j)$ \Comment{Example: $h_{j+4}h_{j+1}c_j = h_{j+3}e_j$}
            \State $m_B \leftarrow \text{CanonicalForm}(h_{j+3}e_j)$
            \If{$m_A \in \text{Map}$ and $m_B \in \text{Map}$}
                \State $\text{Vectors} \leftarrow \text{Vectors} \cup \{ \text{Map}(m_A) + \text{Map}(m_B) \}$
            \EndIf
            \State \Comment{... Loop continues for all 39 relation families ...}
        \EndFor
    \EndIf
    \State \Return Vectors
\EndFunction
\end{algorithmic}

\section*{Python Source Code}

\begin{lstlisting}[language=Python, caption={{\bf SageMath code for calculating Ext group bases}}]
from collections import Counter
from math import log
from sage.all import GF, VectorSpace, matrix
# ===================================================================
# HELPER FUNCTIONS
# ===================================================================
def format_gen(gen_tuple):
    """Formats a generator tuple into a readable algebraic string."""
    type_order = {'P1h1':-2, 'P1h2':-1, 'h':0,'c':1,'d':2,'e':3,'f':4,'g':5,'p':6,'D3':7,'p_prime':8,'n':13,'x':14,'D1':15,'H1':16,'Q3':17,'K':18,'J':19,'T':20,'V':21,'V_prime':22,'U':23}
    counts = Counter(gen_tuple)
    parts = []
    sorted_items = sorted(counts.items(), key=lambda item: (type_order.get(item[0][0], 99), item[0][1]))
    for (gen_type, index), count in sorted_items:
        if gen_type == 'P1h1': type_str = "P^1h_1"
        elif gen_type == 'P1h2': type_str = "P^1h_2"
        else:
            type_str = f"{gen_type}({index})" if gen_type.isupper() else f"{gen_type}_{index}"
            if gen_type in ['p_prime', 'V_prime']: type_str = f"{gen_type.replace('_prime', '')}'_{index}"
        if count > 1: parts.append(f"{type_str}^{count}")
        else: parts.append(type_str)
    return " ".join(parts)

def get_generator_degrees():
    """Returns a dictionary of all generator degrees up to k=5."""
    return {'h':lambda i:2**i,'c':lambda t:2**(t+3)+2**(t+1)+2**t,'d':lambda t:2**(t+4)+2**(t+1),'e':lambda t:2**(t+4)+2**(t+2)+2**t,'f':lambda t:2**(t+4)+2**(t+2)+2**(t+1),'g':lambda t:2**((t-1)+4)+2**((t-1)+3)if t>0 else 0,'p':lambda t:2**(t+5)+2**(t+2)+2**t,'D3':lambda t:2**(t+6)+2**t,'p_prime':lambda t:2**(t+6)+2**(t+3)+2**t,'n':lambda i:2**(i+5)+2**(i+2),'x':lambda i:2**(i+5)+2**(i+3)+2**(i+1),'D1':lambda i:2**(i+5)+2**(i+4)+2**(i+2)+2**i,'H1':lambda i:2**(i+6)+2**(i+1)+2**i,'Q3':lambda i:2**(i+6)+2**(i+3),'K':lambda i:2**(i+7)+2**(i+1),'J':lambda i:2**(i+7)+2**(i+2)+2**i,'T':lambda i:2**(i+7)+2**(i+4)+2**(i+1),'V':lambda i:2**(i+7)+2**(i+5)+2**i,'V_prime':lambda i:2**(i+8)+2**i,'U':lambda i:2**(i+8)+2**(i+3)+2**i,'P1h1':lambda:14,'P1h2':lambda:16}

def canonical_form(gen_tuple):
    """Sorts a generator tuple into a canonical form for consistent identification."""
    type_order = {'P1h1':-2, 'P1h2':-1, 'h':0,'c':1,'d':2,'e':3,'f':4,'g':5,'p':6,'D3':7,'p_prime':8,'n':13,'x':14,'D1':15,'H1':16,'Q3':17,'K':18,'J':19,'T':20,'V':21,'V_prime':22,'U':23}
    return tuple(sorted(list(gen_tuple), key=lambda item: (type_order.get(item[0], 99), item[1])))

def _get_relations(k, potential_gens, gen_map, V, get_vec, max_idx):
    """Helper function to find the full set of Adams relations for a given k."""
    relations = []
    # General h-relations
    if k >= 3:
        for gen in potential_gens:
            h_indices = [p[1] for p in gen if p[0] == 'h']
            for i in range(1, max_idx):
                if h_indices.count(i) >= 3:
                    other_parts = [p for p in gen if p[0] != 'h' or p[1] != i]
                    other_parts.extend([('h',i)] * (h_indices.count(i) % 3))
                    replacement_parts = list(other_parts)
                    for _ in range(h_indices.count(i) // 3):
                        replacement_parts.extend([('h',i-1), ('h',i-1), ('h',i+1)])
                    relations.append(get_vec(gen) + get_vec(canonical_form(tuple(replacement_parts))))
            for i in range(max_idx - 1):
                if i in h_indices and i + 1 in h_indices: relations.append(get_vec(gen))
            for i in range(max_idx - 2):
                if h_indices.count(i) >= 1 and h_indices.count(i + 2) >= 2: relations.append(get_vec(gen))
            for i in range(max_idx - 3):
                if h_indices.count(i) >= 2 and h_indices.count(i + 3) >= 2: relations.append(get_vec(gen))
    # Relations for k>=4 products involving c_i
    if k >= 4:
        for gen in potential_gens:
            c_indices = [p[1] for p in gen if p[0] == 'c']; h_indices = [p[1] for p in gen if p[0] == 'h']
            for c_idx in c_indices:
                for v_rel in [c_idx-1, c_idx, c_idx+2, c_idx+3]:
                    if v_rel in h_indices: relations.append(get_vec(gen))   
    if k == 5:
        for j in range(max_idx):
            relations.append(get_vec(canonical_form((('h',j+4),('h',j+4),('c',j))))); relations.append(get_vec(canonical_form((('h',j+3),('h',j),('c',j+2))))); relations.append(get_vec(canonical_form((('h',j+1),('h',j+1),('c',j)))))
            relations.append(get_vec(canonical_form((('h',j),('d',j+1))))); relations.append(get_vec(canonical_form((('h',j+3),('d',j))))); relations.append(get_vec(canonical_form((('h',j+4),('d',j)))))
            relations.append(get_vec(canonical_form((('h',j),('e',j+1))))); relations.append(get_vec(canonical_form((('h',j+4),('e',j)))))
            relations.append(get_vec(canonical_form((('h',j+1),('f',j))))); relations.append(get_vec(canonical_form((('h',j+3),('f',j))))); relations.append(get_vec(canonical_form((('h',j+4),('f',j)))))
            if j+1 > 0: relations.append(get_vec(canonical_form((('h',j+3),('g',j+1)))))
            relations.append(get_vec(canonical_form((('h',j),('p',j+1))))); relations.append(get_vec(canonical_form((('h',j+1),('p',j))))); relations.append(get_vec(canonical_form((('h',j+2),('p',j))))); relations.append(get_vec(canonical_form((('h',j+4),('p',j))))); relations.append(get_vec(canonical_form((('h',j+5),('p',j)))))
            relations.append(get_vec(canonical_form((('h',j),('D3',j+1))))); relations.append(get_vec(canonical_form((('h',j),('D3',j))))); relations.append(get_vec(canonical_form((('h',j+5),('D3',j))))); relations.append(get_vec(canonical_form((('h',j+6),('D3',j)))))
            relations.append(get_vec(canonical_form((('h',j),('p_prime',j+1))))); relations.append(get_vec(canonical_form((('h',j+2),('p_prime',j))))); relations.append(get_vec(canonical_form((('h',j+3),('p_prime',j))))); relations.append(get_vec(canonical_form((('h',j+6),('p_prime',j)))))
            relations.append(get_vec(canonical_form((('h',j+4),('h',j+1),('c',j)))) + get_vec(canonical_form((('h',j+3),('e',j)))))
            relations.append(get_vec(canonical_form((('h',j+4),('h',j),('c',j+3)))) + get_vec(canonical_form((('h',j+5),('p_prime',j)))))
            relations.append(get_vec(canonical_form((('h',j+5),('h',j+5),('c',j)))) + get_vec(canonical_form((('h',j+1),('p_prime',j)))))
            relations.append(get_vec(canonical_form((('h',j),('d',j+2)))) + get_vec(canonical_form((('h',j+3),('D3',j)))))
            relations.append(get_vec(canonical_form((('h',j+1),('d',j+1)))) + get_vec(canonical_form((('h',j),('p',j)))))
            if j+1 > 0: relations.append(get_vec(canonical_form((('h',j+2),('d',j+1)))) + get_vec(canonical_form((('h',j+4),('g',j+1)))))
            relations.append(get_vec(canonical_form((('h',j+2),('d',j)))) + get_vec(canonical_form((('h',j),('e',j)))))
            relations.append(get_vec(canonical_form((('h',j+1),('e',j)))) + get_vec(canonical_form((('h',j),('f',j)))))
            if j >= 1: relations.append(get_vec(canonical_form((('h',j+1),('e',j)))) + get_vec(canonical_form((('h',j-1),('h',j-1),('c',j+1)))))
            if j+1 > 0: relations.append(get_vec(canonical_form((('h',j+2),('e',j)))) + get_vec(canonical_form((('h',j),('g',j+1)))))
            relations.append(get_vec(canonical_form((('h',j),('f',j+2)))) + get_vec(canonical_form((('h',j+4),('p_prime',j)))))
            relations.append(get_vec(canonical_form((('h',j),('f',j+1)))) + get_vec(canonical_form((('h',j+3),('p',j)))))
            if j+1 > 0: relations.append(get_vec(canonical_form((('h',j+2),('f',j)))) + get_vec(canonical_form((('h',j+1),('g',j+1)))))
            if j+2 > 0 and j+1 > 0: relations.append(get_vec(canonical_form((('h',j+3),('g',j+2)))) + get_vec(canonical_form((('h',j+5),('g',j+1)))))
    return relations

# ===================================================================
# MAIN CALCULATION FUNCTION
# ===================================================================
def calculate_ext_basis(k, n, verbose=True):
    total_degree = k + n
    if verbose: print(f"--- Calculating basis for Ext_A^({k}, {total_degree}) ---")
    
    gen_degrees = get_generator_degrees(); potential_gens = []
    max_idx = int(log(total_degree, 2)) + 7 if total_degree > 1 else 7
    
    if k==1:
        for i in range(max_idx):
            if gen_degrees['h'](i)==total_degree: potential_gens.append((('h',i),))
    elif k==2:
        for i in range(max_idx):
         for j in range(i,max_idx):
          if gen_degrees['h'](i)+gen_degrees['h'](j)==total_degree: potential_gens.append((('h',i),('h',j)))
    elif k==3:
        for i in range(max_idx):
         for j in range(i,max_idx):
          for l in range(j,max_idx):
           if gen_degrees['h'](i)+gen_degrees['h'](j)+gen_degrees['h'](l)==total_degree: potential_gens.append((('h',i),('h',j),('h',l)))
        for t in range(max_idx):
            try:
                if gen_degrees['c'](t)==total_degree: potential_gens.append((('c',t),))
            except (ValueError, IndexError): continue
    elif k==4:
        for i in range(max_idx):
         for j in range(i,max_idx):
          for l in range(j,max_idx):
           for m in range(l,max_idx):
            if gen_degrees['h'](i)+gen_degrees['h'](j)+gen_degrees['h'](l)+gen_degrees['h'](m)==total_degree: potential_gens.append((('h',i),('h',j),('h',l),('h',m)))
        for u in range(max_idx):
         for v in range(max_idx):
          try:
           if gen_degrees['h'](u)+gen_degrees['c'](v)==total_degree: potential_gens.append(canonical_form((('h',u),('c',v))))
          except (ValueError, IndexError): continue
        for gen_type in ['d','e','f','g','p','D3','p_prime']:
            for t_idx in range(max_idx):
                if gen_type=='g' and t_idx==0: continue
                try:
                 if gen_degrees[gen_type](t_idx)==total_degree: potential_gens.append(((gen_type,t_idx),))
                except (ValueError, IndexError): continue
    elif k==5:
        for i1 in range(max_idx):
         for i2 in range(i1,max_idx):
          for i3 in range(i2,max_idx):
           for i4 in range(i3,max_idx):
            for i5 in range(i4,max_idx):
             if sum(gen_degrees['h'](i) for i in [i1,i2,i3,i4,i5])==total_degree: potential_gens.append(canonical_form((('h',i1),('h',i2),('h',i3),('h',i4),('h',i5))))
        for t in range(max_idx):
         for i in range(max_idx):
          for j in range(i,max_idx):
           if gen_degrees['c'](t)+gen_degrees['h'](i)+gen_degrees['h'](j)==total_degree: potential_gens.append(canonical_form((('c',t),('h',i),('h',j))))
        for i in range(max_idx):
         for j in range(i,max_idx):
          for t in range(max_idx):
           if gen_degrees['h'](i)+gen_degrees['h'](j)+gen_degrees['c'](t)==total_degree: potential_gens.append(canonical_form((('h',i),('h',j),('c',t))))
        for gen_type in ['d','e','f','g','p','D3','p_prime']:
            for i in range(max_idx):
                for t in range(max_idx):
                    if gen_type=='g' and t==0: continue
                    if gen_degrees['h'](i)+gen_degrees[gen_type](t)==total_degree: potential_gens.append(canonical_form((('h',i),(gen_type,t))))
        k5_indecomp_types = ['n','x','D1','H1','Q3','K','J','T','V','V_prime','U']
        for gen_type in k5_indecomp_types:
            for t_idx in range(max_idx):
                if gen_degrees[gen_type](t_idx)==total_degree: potential_gens.append(((gen_type,t_idx),))
        if total_degree == gen_degrees['P1h1'](): potential_gens.append((('P1h1',0),))
        if total_degree == gen_degrees['P1h2'](): potential_gens.append((('P1h2',0),))
    else:
        print(f"Calculation for k > 5 is not supported.")
        return

    potential_gens = sorted(list(set(potential_gens)))
    if not potential_gens:
        if verbose and n>0: print("No potential generators found. Ext group is trivial.")
        return
        
    V = VectorSpace(GF(2), len(potential_gens)); gen_map = {gen: i for i, gen in enumerate(potential_gens)}
    if verbose: print(f"\nFound {V.dimension()} potential generators (before relations):"); [print(f"  {format_gen(g)}") for g in potential_gens]
    
    def get_vec(gen_tuple):
        canon_gen = canonical_form(gen_tuple)
        if canon_gen in gen_map:
            idx = gen_map[canon_gen]; coeffs = [0] * V.dimension(); coeffs[idx] = 1; return V(coeffs)
        return V.zero()
    
    relations = _get_relations(k, potential_gens, gen_map, V, get_vec, max_idx)

    if not relations: relation_matrix = matrix(GF(2), 0, V.dimension())
    else:
        non_zero_relations = [r for r in relations if not r.is_zero()]
        if not non_zero_relations: relation_matrix = matrix(GF(2), 0, V.dimension())
        else: relation_matrix = matrix(GF(2), non_zero_relations)

    kernel = relation_matrix.transpose().kernel(); final_basis_dim = kernel.dimension()
    
    if verbose:
        print(f"\nDimension of relation space: {relation_matrix.rank()}"); print(f"Dimension of Ext_A^({k}, {total_degree}) = {final_basis_dim}")
        
        relation_space = relation_matrix.row_space()
        echelon_matrix = relation_matrix.echelon_form()
        if not echelon_matrix.is_zero():
            print("\nSimplified Relations:")
            for rel_vec in echelon_matrix:
                if rel_vec.is_zero(): continue
                indices = rel_vec.nonzero_positions()
                if len(indices) == 1: lhs = format_gen(potential_gens[indices[0]]); print(f"  -> {lhs} = 0")
                else: lhs = format_gen(potential_gens[indices[0]]); rhs = " + ".join([format_gen(potential_gens[i]) for i in indices[1:]]); print(f"  -> {lhs} = {rhs}")

        if final_basis_dim > 0:
            print("\nBasis elements:")
            basis_reps = [potential_gens[v.nonzero_positions()[0]] for v in kernel.basis()]
            for b_rep in sorted(basis_reps):
                v_b = get_vec(b_rep); equivalents = [b_rep]
                for g in potential_gens:
                    if g == b_rep: continue
                    v_g = get_vec(g)
                    if (v_b - v_g) in relation_space: equivalents.append(g)
                sorted_equivalents = sorted(equivalents, key=lambda gen: (format_gen(gen).find('='), format_gen(gen)))
                print(f"{ ' = '.join([format_gen(g) for g in sorted_equivalents])}")

\end{lstlisting}

\medskip

{\bf Example:} To verify certain results from Nguyen Sum's paper~\cite{Sum} for degrees $n = 2^{s+1} - m$ with $m = 2,\, 3$, as well as results from our previous works~\cite{Phuc1, Phuc3} for the degree $n = 2^{s+t+u} + 2^{s+t} + 2^{s} - 3$, and to provide sample computations for $k = 5$, we define the following function along with example calculations for rank $k = 4$ at the above degrees and for rank $k = 5$ with $n \in \big\{62,\, 128,\, 256\big\}.$ This function should be inserted immediately after the line of code \texttt{print(f"{ ' = '.join([format\_gen(g) for g in sorted\_equivalents])}")} in the \textsc{SageMath} script above.

\medskip

\begin{lstlisting}[language=Python]
def calculate_for_n_s(k, s, m):
    print(f"====================================================")
    print(f"Calculating for Case n = 2^{s+1} - m with k={k}, s={s}, m={m}")
    print(f"====================================================\n")
    if not (s >= 0 and m in [2, 3]):
        print("Warning: s should be >= 0 and m must be 2 or 3."); return
    if not (1 <= k <= 5):
        print(f"Warning: This script is only configured for k between 1 and 5. Received k={k}."); return
    n = 2**(s+1) - m
    if n < 0:
        print(f"--- Skipping case because n = {n} is negative ---\n"); return
    print(f"----Case: n = 2^{s+1} - {m} = {n} ----")
    calculate_ext_basis(k=k, n=n, verbose=True); print("\n")

def calculate_for_n_stu(k, s, t, u):
    print(f"====================================================")
    print(f"Calculating for Case n_stu with k={k}, (s={s}, t={t}, u={u})")
    print(f"====================================================\n")
    if not (s >= 0 and t >= 0 and u >= 0):
        print("Warning: s, t, u should be >= 0."); return
    if not (1 <= k <= 5):
        print(f"Warning: This script is only configured for k between 1 and 5. Received k={k}."); return
    n = 2**(s + t + u) + 2**(s + t) + 2**s - 3
    print(f"--Case: n = 2^({s+t+u}) + 2^({s+t}) + 2^{s} - 3 = {n} --")
    calculate_ext_basis(k=k, n=n, verbose=True); print("\n")

print("--- Running calculations for n = 2^{s+1} - m ---")
calculate_for_n_s(k=4, s=6, m=2)
calculate_for_n_s(k=4, s=5, m=3)
print("\n--- Running calculation for n_stu ---")
calculate_for_n_stu(k=4, s=2, t=1, u=2)
print("--- Calculation for k=5, n=62 ---")
calculate_ext_basis(k=5, n=62)
print("--- Calculation for k=5, n=128 ---")
calculate_ext_basis(k=5, n=128)
print("--- Calculation for k=5, n=256 ---")
calculate_ext_basis(k=5, n=256)
\end{lstlisting}

\medskip

The algorithm yields the following results for the cases in the "\texttt{SAMPLE CALCULATION}" section:

\begin{lstlisting}
--- Running calculations for n = 2^{s+1} - m ---
========================================================================
Calculating for Case n = 2^7 - m with k=4, s=6, m=2
========================================================================
-------- Case: n = 2^7 - 2 = 126 --------
--- Calculating basis for Ext_A^(4, 130) ---
Found 3 potential generators (before relations):
  D_3(1)
  h_0^2h_6^2
  h_1h_5^2h_6
Dimension of relation space: 1
Dimension of Ext_A^(4, 130) = 2
Simplified Relations:
  -> h_1h_5^2h_6 = 0
Basis elements:
  D_3(1)
  h_0^2h_6^2
========================================================================
Calculating for Case n = 2^6 - m with k=4, s=5, m=3 
========================================================================
-------- Case: n = 2^6 - 3 = 61 --------
--- Calculating basis for Ext_A^(4, 65) ---
Found 2 potential generators (before relations):
  D_3(0)
  h_0h_4^2h_5
Dimension of relation space: 1
Dimension of Ext_A^(4, 65) = 1
Simplified Relations:
  -> h_0h_4^2h_5 = 0
Basis elements:
  D_3(0)

--- Running calculation for n_stu ---
========================================================================
Calculating for Case n_stu with k=4, (s=2, t=1, u=2)   
========================================================================
-------- Case: n = 2^(5) + 2^(3) + 2^2 - 3 = 41 --------
--- Calculating basis for Ext_A^(4, 45) ---
Found 2 potential generators (before relations):
  h_0c_2
  h_0h_2h_3h_5
Dimension of relation space: 1
Dimension of Ext_A^(4, 45) = 1
Simplified Relations:
  -> h_0h_2h_3h_5 = 0
Basis elements:
  h_0c_2

--- Calculation for k=5, n=62 ---
--- Calculating basis for Ext_A^(5, 67) ---
Found 4 potential generators (before relations):
  H_1(0)
  h_0^3h_5^2
  h_0h_1h_4^2h_5
  h_1D_3(0)

Dimension of relation space: 1
Dimension of Ext_A^(5, 67) = 3
Simplified Relations:
  -> h_0h_1h_4^2h_5 = 0

Basis elements:
  H_1(0)
  h_0^3h_5^2
  h_1D_3(0)

--- Calculation for k=5, n=128 ---
--- Calculating basis for Ext_A^(5, 133) ---

Found 4 potential generators (before relations):
  J(0)
  h_0^3h_1h_7
  h_0h_1^2h_6^2
  h_0h_2h_5^2h_6

Dimension of relation space: 3
Dimension of Ext_A^(5, 133) = 1
Simplified Relations:
  -> h_0^3h_1h_7 = 0
  -> h_0h_1^2h_6^2 = 0
  -> h_0h_2h_5^2h_6 = 0

Basis elements:
  J(0)

--- Calculation for k=5, n=256 ---
--- Calculating basis for Ext_A^(5, 261) ---

Found 4 potential generators (before relations):
  h_0D_3(2)
  h_0^3h_1h_8
  h_0 h_1^2h_7^2
  h_0 h_2 h_6^2h_7

Dimension of relation space: 3
Dimension of Ext_A^(5, 261) = 1

Simplified Relations:
  -> h_0^3h_1h_8 = 0
  -> h_0h_1^2h_7^2 = 0
  -> h_0h_2h_6^2h_7 = 0

Basis elements:
  h_0D_3(2)
\end{lstlisting}

\medskip
\newpage
Thus, from the above algorithmic output, we obtain the following results:

\medskip

$\bullet$ For $k  =4$ and degree $ n  =2^{6+1}-2 = 126,$ we have

\medskip

\begin{tcolorbox}[colback=white!5!white, colframe=black!75!black, title={Computation of $ Ext_{\mathcal A}^{4, 4+126}(\mathbb Z/2, \mathbb Z/2)$}]
$$ Ext_{\mathcal A}^{4, 4+126}(\mathbb Z/2, \mathbb Z/2) = \langle D_3(1), h_0^{2}h_6^{2} \rangle.$$
\end{tcolorbox}

\medskip

$\bullet$ For $k  =4$ and degree $ n  =2^{5+1}-3 = 61$ we have

\medskip

\begin{tcolorbox}[colback=white!5!white, colframe=black!75!black, title={Computation of $Ext_{\mathcal A}^{4, 4+61}(\mathbb Z/2, \mathbb Z/2) $}]
$$ Ext_{\mathcal A}^{4, 4+61}(\mathbb Z/2, \mathbb Z/2) = \langle D_3(0) \rangle.$$
\end{tcolorbox}

These output results from the algorithm indicate that the computations in Nguyen Sum's original paper~\cite{Sum} are incorrect for the degrees $n = 2^{s+1} - 2$ with $s = 6$, and $n = 2^{s+1} - 3$ with $s = 5.$ These results were later revised by Nguyen Sum in the arXiv version:1710.07895.

\medskip

$\bullet$ For $k  =4$ and degree $n_{2,\, 1,\, 2} = 41,$ we have 

\medskip

\begin{tcolorbox}[colback=white!5!white, colframe=black!75!black, title={Computation of $Ext_{\mathcal A}^{4, 4+n_{2,\, 1,\, 2}}(\mathbb Z/2, \mathbb Z/2) $}]
$$ Ext_{\mathcal A}^{4, 4+n_{2,\, 1,\, 2}}(\mathbb Z/2, \mathbb Z/2) = \langle h_0c_2 \rangle.$$
\end{tcolorbox}
\medskip

In a similar manner, the reader may substitute any values of $(s, t, u)$ into the code above to obtain the corresponding results, matching those we previously computed by hand in~\cite{Phuc1, Phuc3}. More precisely, 

\medskip

\begin{tcolorbox}[colback=white!5!white, colframe=black!75!black, title={Computation of ${Ext}^{4, 4+n_{s,\, t,\, u}}_{\mathcal A}(\mathbb Z/2, \mathbb Z/2) $}]
$$
{Ext}^{4, 4+n_{s,\, t,\, u}}_{\mathcal A}(\mathbb Z/2, \mathbb Z/2) = \left\{\begin{array}{ll}
\langle h_0c_s \rangle&\mbox{if $s \geq 2$, $t = 1$ and $u  =2$},\\
\langle h_{u+3}c_0 \rangle&\mbox{if $s = 1$, $t =  2$ and $u\geq 1$},\\
\langle h_0h_sh_{s+t}h_{s+t+u} \rangle&\mbox{if $s \geq 2$, $t\geq 2$ and $u\geq 2$},\\
0&\mbox{otherwise}.
\end{array}\right.
$$
\end{tcolorbox}

Note that the remaining cases of $(s,\, t,\, u)$ for which ${Ext}^{4,\, 4+n_{s,\, t,\, u}}_{\mathcal A}(\mathbb{Z}/2, \mathbb{Z}/2) = 0$ are of the following forms:

\medskip

$-$ $s = 1$, $t \geq 1,\, t\neq 2$ and $u \geq 1;$ 

$-$ $s = 2$, $t = 1$ and $u \geq 1,\, u\neq 2;$ 

$-$ $s = 2$, $t \geq 2$ and $u = 1;$

$-$ $s \geq 3$, $t = 1$ and $u \geq 1,\, u\neq 2;$ 

$-$ $s \geq 3$, $t \geq 2$ and $u =1.$

\medskip

$\bullet$ For $k  =5$ and degree $n = 62,$ we have

\medskip

\begin{tcolorbox}[colback=white!5!white, colframe=black!75!black, title={Computation of $Ext_{\mathcal A}^{5, 5+62}(\mathbb Z/2, \mathbb Z/2) $}]
$$  Ext_{\mathcal A}^{5, 5+62}(\mathbb Z/2, \mathbb Z/2) = \langle  \{H_1(0),\, h_0^3 h_5^2,\, h_1 D_3(0)\} \rangle.$$
\end{tcolorbox}

\medskip
\newpage
$\bullet$ For $k  =5$ and degree $n = 128,$ we have

\medskip

\begin{tcolorbox}[colback=white!5!white, colframe=black!75!black, title={Computation of $Ext_{\mathcal A}^{5, 5+128}(\mathbb Z/2, \mathbb Z/2) $}]
$$  Ext_{\mathcal A}^{5, 5+128}(\mathbb Z/2, \mathbb Z/2) = \langle  J_0 \rangle.$$
\end{tcolorbox}

\medskip

$\bullet$ For $k  =5$ and degree $n = 256,$ we have

\medskip

\begin{tcolorbox}[colback=white!5!white, colframe=black!75!black, title={Computation of $Ext_{\mathcal A}^{5, 5+256}(\mathbb Z/2, \mathbb Z/2) $}]
$$  Ext_{\mathcal A}^{5, 5+256}(\mathbb Z/2, \mathbb Z/2) = \langle  h_0D_3(2) \rangle.$$
\end{tcolorbox}

\medskip

\section{Appendix: Code for computing the generalized results of the generators of ${\rm Ext}^{4,\, 4+n_{s,\, t,\, u}}_{\mathcal{A}}(\mathbb{Z}/2, \mathbb{Z}/2)$}

\begin{lstlisting}[language=Python, caption={{\bf SageMath code for computing the generators of ${\rm Ext}^{4,\, 4+n_{s,\, t,\, u}}_{\mathcal{A}}(\mathbb{Z}/2, \mathbb{Z}/2)$.}}]
from collections import Counter
from math import log
from sage.all import GF, VectorSpace, matrix
import re
# ===================================================================
# HELPER FUNCTIONS
# ===================================================================
def format_gen(gen_tuple):
    """Formats a generator tuple into a readable algebraic string."""
    type_order = {'h': 0, 'c': 1}
    counts = Counter(gen_tuple)
    parts = []
    sorted_items = sorted(counts.items(), key=lambda item: (type_order.get(item[0][0], 99), item[0][1]))
    for (gen_type, index), count in sorted_items:
        type_str = f"{gen_type}_{index}"
        if count > 1:
            parts.append(f"{type_str}^{count}")
        else:
            parts.append(type_str)
    return " ".join(parts)

def get_generator_degrees():
    """Returns a dictionary of generator degrees for h and c."""
    return {'h': lambda i: 2**i, 'c': lambda t: 2**(t + 3) + 2**(t + 1) + 2**t}

def canonical_form(gen_tuple):
    """Sorts a generator tuple into a canonical form for consistent identification."""
    type_order = {'h': 0, 'c': 1}
    return tuple(sorted(list(gen_tuple), key=lambda item: (type_order.get(item[0], 99), item[1])))

def calculate_ext_basis(k, n, verbose=False):
    """
    Calculates the Z/2-basis for Ext_A^(k, k+n) by applying Adams relations.
    """
    total_degree = k + n
    gen_degrees = get_generator_degrees()
    potential_gens = []
    max_idx = int(log(total_degree, 2)) + 7 if total_degree > 1 else 7

    if k == 4:
        # 4h generators
        for i in range(max_idx):
            for j in range(i, max_idx):
                for l in range(j, max_idx):
                    for m in range(l, max_idx):
                        if gen_degrees['h'](i) + gen_degrees['h'](j) + gen_degrees['h'](l) + gen_degrees['h'](m) == total_degree:
                            potential_gens.append((('h', i), ('h', j), ('h', l), ('h', m)))

        # h+c generators
        for u in range(max_idx):
            for v in range(max_idx):
                try:
                    if gen_degrees['h'](u) + gen_degrees['c'](v) == total_degree:
                        potential_gens.append(canonical_form((('h', u), ('c', v))))
                except (ValueError, IndexError):
                    continue

    potential_gens = sorted(list(set(potential_gens)))
    if not potential_gens:
        return 0, []

    V = VectorSpace(GF(2), len(potential_gens))
    gen_map = {gen: i for i, gen in enumerate(potential_gens)}

    def get_vec(gen_tuple):
        canon_gen = canonical_form(gen_tuple)
        if canon_gen in gen_map:
            idx = gen_map[canon_gen]
            coeffs = [0] * V.dimension()
            coeffs[idx] = 1
            return V(coeffs)
        return V.zero()

    #  Relations for k=4
    relations = []
    for gen in potential_gens:
        h_indices = [p[1] for p in gen if p[0] == 'h']
        c_indices = [p[1] for p in gen if p[0] == 'c']

        # h-generator relations
        for i in range(max_idx - 1):
            if i in h_indices and i + 1 in h_indices: relations.append(get_vec(gen)) # h_i * h_{i+1} = 0
        for i in range(max_idx):
            if h_indices.count(i) >= 3: relations.append(get_vec(gen))
        for i in range(max_idx - 2):
            if h_indices.count(i) >= 1 and h_indices.count(i + 2) >= 2: relations.append(get_vec(gen)) # h_i * h_{i+2}^2 = 0
        for i in range(max_idx - 3):
            if h_indices.count(i) >= 2 and h_indices.count(i + 3) >= 2: relations.append(get_vec(gen)) # h_i^2 * h_{i+3}^2 = 0

        # h-c relations
        if len(h_indices) == 1 and len(c_indices) == 1:
            h_idx, c_idx = h_indices[0], c_indices[0]
            if h_idx in [c_idx - 1, c_idx, c_idx + 2, c_idx + 3]: relations.append(get_vec(gen))

    if not relations:
        relation_matrix = matrix(GF(2), 0, V.dimension())
    else:
        non_zero_relations = [r for r in relations if not r.is_zero()]
        if not non_zero_relations:
            relation_matrix = matrix(GF(2), 0, V.dimension())
        else:
            relation_matrix = matrix(GF(2), non_zero_relations)

    kernel = relation_matrix.transpose().kernel()
    final_basis_dim = kernel.dimension()

    return final_basis_dim, kernel.basis() if final_basis_dim > 0 else []

def extract_numerical_data(generator_str):
    """Extracts numerical indices from a formatted generator string."""
    h_matches = re.findall(r'h_(\d+)', generator_str)
    c_matches = re.findall(r'c_(\d+)', generator_str)
    return [int(x) for x in h_matches], [int(x) for x in c_matches]

def find_canonical_description(value, s, t, u):
    """Finds the most canonical parametric description for a value."""
    if value == 0: return "0"
    if value == s: return "s"
    if value == t: return "t"
    if value == u: return "u"
    
    for const in range(1, 6):
        if value == u + const: return f"u+{const}"
        if value == t + const: return f"t+{const}"
        if value == s + const: return f"s+{const}"
    
    if value == s + t: return "s+t"
    if value == s + u: return "s+u"
    if value == t + u: return "t+u"
    if value == s + t + u: return "s+t+u"
    
    return str(value)

def analyze_conditions_optimally(conditions):
    """Analyzes parameter conditions to produce clean mathematical descriptions."""
    if not conditions: return "no conditions"
    
    s_values = sorted(set(cond[0] for cond in conditions))
    t_values = sorted(set(cond[1] for cond in conditions))
    u_values = sorted(set(cond[2] for cond in conditions))
    
    def analyze_parameter_values(values, param_name):
        if len(values) == 1: return f"{param_name} = {values[0]}"
        elif len(values) >= 2: return f"{param_name} >= {min(values)}" 
        return ""
    
    s_desc = analyze_parameter_values(s_values, 's')
    t_desc = analyze_parameter_values(t_values, 't')
    u_desc = analyze_parameter_values(u_values, 'u')
    
    return ", ".join(filter(None, [s_desc, t_desc, u_desc]))

def evaluate_parametric_expression(expr, s, t, u):
    """Evaluates a parametric expression string with given s, t, u values."""
    try:
        return eval(expr, {"s": s, "t": t, "u": u})
    except Exception:
        return None

def discover_pure_patterns(nonzero_cases):
    """
    Analyzes computational results to discover and consolidate general patterns.
    """
    index_pattern_groups = {}
    
    for case in nonzero_cases:
        s, t, u = case['s'], case['t'], case['u']
        gen_str = case['generator_pattern']
        h_indices, c_indices = extract_numerical_data(gen_str)
        
        pattern_key = "unknown"
        if len(h_indices) == 1 and len(c_indices) == 1:
            h_desc = find_canonical_description(h_indices[0], s, t, u)
            c_desc = find_canonical_description(c_indices[0], s, t, u)
            pattern_key = f"h_{{{h_desc}}}c_{{{c_desc}}}"
        elif len(h_indices) == 4 and len(c_indices) == 0:
            h_sorted = sorted(h_indices)
            if h_sorted == sorted([0, s, s + t, s + t + u]):
                pattern_key = "h_0h_sh_{s+t}h_{s+t+u}"
            else:
                pattern_key = "h_" + "_".join([str(idx) for idx in h_sorted])
        else:
            continue
        
        if pattern_key not in index_pattern_groups:
            index_pattern_groups[pattern_key] = []
        index_pattern_groups[pattern_key].append(case)
    
    final_patterns = []
    for pattern, cases in index_pattern_groups.items():
        if len(cases) < 2:
            continue
        
        is_valid = True
        for case in cases:
            s, t, u = case['s'], case['t'], case['u']
            h_indices, c_indices = extract_numerical_data(case['generator_pattern'])
            
            if 'h_{' in pattern and '}c_{' in pattern:
                match = re.match(r'h_\{([^}]+)\}c_\{([^}]+)\}', pattern)
                if match and len(h_indices) == 1 and len(c_indices) == 1:
                    h_expr, c_expr = match.groups()
                    expected_h = evaluate_parametric_expression(h_expr, s, t, u)
                    expected_c = evaluate_parametric_expression(c_expr, s, t, u)
                    if expected_h != h_indices[0] or expected_c != c_indices[0]:
                        is_valid = False
                        break
        
        if is_valid:
            conditions = [(case['s'], case['t'], case['u']) for case in cases]
            final_patterns.append({
                'pattern': pattern,
                'condition': analyze_conditions_optimally(conditions),
                'case_count': len(cases)
            })
    
    final_patterns.sort(key=lambda x: -x['case_count'])
    return final_patterns

def generate_clean_theorem(patterns, total_cases, s_max, t_max, u_max):
    """Generates a formatted theorem-style output from the discovered patterns."""
    
    print("="*80)
    print(f"General Form of Generators of Ext^{{4, 4+n_{{s,t,u}}}}_A")
    print(f"Discovered from {total_cases} computed cases for 1 <= s <= {s_max}, 1 <= t <= {t_max}, 1 <= u <= {u_max}")
    print("="*80)
    
    if not patterns:
        print("No significant patterns discovered.")
        return
    
    print("Ext^{4, 4+n_{s,t,u}}_A(F_2, F_2) = {")
    
    max_pattern_length = max(len(f"<{p['pattern']}>") for p in patterns) if patterns else 0 
    
    for pattern_info in patterns:
        pattern_str = f"<{pattern_info['pattern']}>" 
        padding = max_pattern_length - len(pattern_str) + 4
        print(f"  {pattern_str}{' ' * padding}if {pattern_info['condition']}")
    
    padding_otherwise = max_pattern_length - 1 + 4
    print(f"  0{' ' * padding_otherwise}otherwise")
    print("}")
    print("="*80)

def run_pure_discovery(s_max=10, t_max=10, u_max=10, show_details=False):
    """Main discovery function that computes cases and derives general patterns."""
    
    all_cases = []
    
    for s in range(1, s_max + 1):
        for t in range(1, t_max + 1):
            for u in range(1, u_max + 1):
                n = int(2**(s + t + u) + 2**(s + t) + 2**s - 3)
                dim, basis = calculate_ext_basis(k=4, n=n, verbose=False)
                
                case_info = {'s': s, 't': t, 'u': u, 'n': n, 'dim': dim}
                
                if dim == 1 and basis:
                    total_degree = 4 + n
                    gen_degrees = get_generator_degrees()
                    potential_gens = []
                    max_idx = int(log(total_degree, 2)) + 7 if total_degree > 1 else 7
                    
                    for i in range(max_idx):
                     for j in range(i,max_idx):
                      for l in range(j,max_idx):
                       for m in range(l,max_idx):
                        if gen_degrees['h'](i)+gen_degrees['h'](j)+gen_degrees['h'](l)+gen_degrees['h'](m)==total_degree:
                            potential_gens.append((('h',i),('h',j),('h',l),('h',m)))
                    for v_idx in range(max_idx):
                     for w_idx in range(max_idx):
                      try:
                       if gen_degrees['h'](v_idx)+gen_degrees['c'](w_idx)==total_degree: 
                           potential_gens.append(canonical_form((('h',v_idx),('c',w_idx))))
                      except (ValueError, IndexError): continue
                    potential_gens = sorted(list(set(potential_gens)))
                    
                    if potential_gens and len(basis) > 0:
                        basis_vector = basis[0]
                        for i, coeff in enumerate(basis_vector):
                            if coeff != 0:
                                primary_gen = potential_gens[i]
                                case_info['generator_pattern'] = format_gen(primary_gen)
                                break
                
                all_cases.append(case_info)
    
    nonzero_cases = [case for case in all_cases if case['dim'] > 0 and 'generator_pattern' in case]
    
    if show_details:
        print(f"\nComputed {len(all_cases)} total cases.")
        print(f"Found {len(nonzero_cases)} non-zero, one-dimensional cases.")
        print("\nShowing first 10 non-zero cases:")
        for case in nonzero_cases[:10]:
            print(f"  (s={case['s']}, t={case['t']}, u={case['u']}) -> {case['generator_pattern']}")
        if len(nonzero_cases) > 10:
            print(f"  ... and {len(nonzero_cases) - 10} more.")

    patterns = discover_pure_patterns(nonzero_cases)
    
    generate_clean_theorem(patterns, len(all_cases), s_max, t_max, u_max)
    
    return all_cases, patterns
# ===================================================================
# EXECUTION
# ===================================================================
print("Starting the computation of Ext^{4, 4+n}_A...")
print("Case n = n_{s,t,u} = 2^{s+t+u} + 2^{s+t} + 2^{s} - 3, for s, t, u >= 1")
print()
S_MAX, T_MAX, U_MAX = 10, 10, 10
cases, patterns = run_pure_discovery(s_max=S_MAX, t_max=T_MAX, u_max=U_MAX, show_details=False)
print("\nComputation successfully completed!")
\end{lstlisting}

\medskip

The results produced by the above algorithm coincide with the manually computed results in our previous works~\cite{Phuc1, Phuc3}:

\medskip

\begin{lstlisting}
Starting the computation of Ext^{4, 4+n}_A...
Case n = n_{s,t,u} = 2^{s+t+u} + 2^{s+t} + 2^{s} - 3, for s, t, u >= 1
========================================================================
General Form of Generators of Ext^{4, 4+n_{s,t,u}}_A
Discovered from 1000 computed cases for 
1 <= s <= 10, 1 <= t <= 10, 1 <= u <= 10
========================================================================
Ext^{4, 4+n_{s,t,u}}_A(F_2, F_2) = {
  <h_0h_sh_{s+t}h_{s+t+u}>    if s >= 2, t >= 2, u >= 2
  <h_{u+3}c_{0}>              if s = 1, t = 2, u >= 2
  <h_{0}c_{s}>                if s >= 2, t = 1, u = 2
    0                         otherwise
}
========================================================================
Computation successfully completed!
\end{lstlisting}

\medskip

\section*{Acknowledgments}
The author thanks Professor John Rognes (University of Oslo, Norway) for his helpful email correcting the terminology regarding "Adem relations" of $\mathrm{Ext}_{\mathcal A}^{k, k+*}(\mathbb{Z}/2, \mathbb{Z}/2)$ for $k\leq 3.$

\end{document}